\documentclass{amsart}

\usepackage{amssymb}
\usepackage{amsmath}
\usepackage{hyperref}

\newcommand{\half}{\frac{1}{2}}

\newcommand{\summ}{\mathop{{\sum}^{\star}}}

\numberwithin{equation}{section}

\newtheorem{theorem}{Theorem}[section]

\newtheorem{lemma}[theorem]{Lemma}

\makeatletter
\@namedef{subjclassname@2010}{%
  \textup{2010} Mathematics Subject Classification}
\makeatother

\begin{document}

\title{The divisor function in arithmetic progressions modulo prime powers}

\author{Rizwanur Khan}
\address{
Science Program\\ Texas A\&M University at Qatar\\ Doha, Qatar}
\email{rizwanur.khan@qatar.tamu.edu }


\begin{abstract} We study the average value of the divisor function $\tau(n)$ for $n\le x$ with $n \equiv a \bmod q$. The divisor function is known to be evenly distributed over arithmetic progressions for all $q$ that are a little smaller than $x^{2/3}$. We show how to go past this barrier when $q=p^k$ for odd primes $p$ and any fixed integer $k\ge 7$.

 \end{abstract}

\maketitle

\section{Introduction}

The study of arithmetic functions along arithmetic progressions has a long and extensive history in number theory. Perhaps the most famous example of this problem is counting the number of primes up to $x$ that are congruent to $a \bmod q$ for some $(a,q)=1$. The classical Siegel-Walfisz theorem \cite{wal} gives the expected asymptotic for this number when $q\le (\log x)^N$ for any fixed $N$. Counting primes in arithmetic progressions for larger moduli is connected to deep unsolved problems on the zeros of Dirichlet $L$-functions (see \cite[Chapter 22]{dav}). One lesson to take away is that such problems are very sensitive to the relative sizes of $x$ and $q$.

The average value of the divisor function $\tau(n)$ in arithmetic progressions is another well known classical problem on which many important applications hinge. For example, as M. Young explains in \cite[page 4]{you}, this problem is central to the fourth moment of Dirichlet $L$-functions. Let $(a,q)=1$ and $x>q^{\theta}$. It is conjectured that if $\theta> 1$, then
\begin{align}
\label{conj} \sum_{\substack{n\le x\\ n\equiv a \bmod q}} \tau(n) - \frac{1}{\phi(q)} \sum_{\substack{n\le x\\ (n,q)=1}} \tau(n) \ll \frac{x^{1-\delta}}{q}
\end{align}
for some $\delta>0$ depending only on $\theta$, where $\phi(q)$ is the Euler totient function. In other words, if the modulus $q$ is not too large, it is expected that the divisor function is evenly distributed over residue classes mod $q$. Consider the left hand side of (\ref{conj}). As we will see below, if one picks out the residue class mod $q$ using additive characters and applies Poisson summation, then one arrives at a certain sum of Kloosterman sums. At this point, applying Weil's bound for each individual Kloosterman sum yields an admissible error term for $\theta>\frac{3}{2}$. In fact a smaller error term is possible; see \cite{ponvau}. To do better, one must seek cancellation between Kloosterman sums instead of bounding them absolutely. This seems to be a very difficult problem and the conjecture (\ref{conj}) remains unsolved for every value of $\theta\le \frac{3}{2}$.

Several authors have provided evidence for (\ref{conj}) beyond the barrier of $\theta=\frac{3}{2}$ by showing that it  holds in an average sense, where the averaging is performed over $q$ or $a$. See \cite{blo, ban, fouiwa, fou}. Of these, in \cite{fouiwa}, Fouvry and Iwaniec work only with moduli having a special factorization. In a recent paper, Irving \cite{irv} was the first to consider some individual moduli beyond the range given by Weil's bound. He showed that (\ref{conj}) holds for $x>q^{1.49}$, provided that $q$ is square-free and has only small prime factors. On page 6679 of his paper, Irving stressed the importance to his method of having $q$ square-free.

In this paper we consider prime power moduli. While this covers an important case which Irving could not treat, the real motivation for this choice of special moduli is its historical significance. A long line of papers \cite{mil, blomil2, blomil, mun, hea, fuj, iwa, gal, bar, pos} concerning Dirichlet characters, $L$-functions, and primes in arithmetic progressions are specialized to, or build upon work specialized to, prime power moduli. We prove the following result.

\begin{theorem}
 \label{main} Let $q=p^k$ for an odd prime $p$ and a fixed integer $k\ge 7$. There exist some constants $\eta>0$ and $\delta>0$, depending on $k$,  such that (\ref{conj}) holds for $x>q^{\frac{3}{2}-\eta}$. The implied constant in (\ref{conj}) depends on $k$.
\end{theorem}
\noindent Thus the main purpose of this paper is show how the barrier of $\theta=\frac{3}{2}$ may be broken for sufficiently powerful moduli $p^k$ as $p\to \infty$. The idea of our method is completely different from that of Irving. As alluded to above, we must obtain cancellation in a sum of Kloosterman sums. To do this we will use the fact that Kloosterman sums to prime power moduli have a special explicit evaluation. Then we will use the theory of exponential sums (more specifically, Weyl differencing) to obtain the required cancellation. As an application, it would interesting to try to use our result to prove an asymptotic with a power saving error term for the fourth moment of Dirichlet $L$-functions to prime power moduli. This would be an alternate method to the one in \cite{blomil}. In an effort to keep the argument transparent we have not computed the values for $\eta$ and $\delta$ in terms of $k$, which would have been minimal (about as small as $2^{-k}$, due to the exponent from Weyl differencing; see (\ref{lastt})).

Theorem \ref{main} will be proven from the following result, which beats Weil's bound for arbitrarily short averages of Kloosterman sums to sufficiently powerful moduli.
\begin{theorem}\label{main2} Let $\lambda>0$ be fixed. If $k> \min\{ \frac{3}{2\lambda},4 \}$ is a fixed integer,  $q=p^k$ for an odd prime $p$,  $N\ge q^{\lambda}$ and $(\beta,p)=1$, then there exists $\delta>0$, depending on $\lambda$ and $k$,  such that
\begin{align}
\label{main2line} \sum_{1\le n \le N} S(n,\beta;q) \ll N q^{\frac{1}{2}-\delta},
\end{align}
where $S(n,\beta;q)$ denotes the Kloosterman sum. The implied constant depends on $\lambda$ and $k$.
\end{theorem}
\noindent This result with $N$ about size $q^\frac{1}{4}$ is needed for Theorem $\ref{main}$.

\section{Proof of Theorem \ref{main2}}

For the terms on the left hand side of (\ref{main2line}) with $p|n$ we have, using Weil's bound $S(n,\beta ;q)\ll q^\half$, that
\begin{align*}
\sum_{\substack{1\le n \le N\\p|n}} S(n,\beta ;q) \ll N q^{\frac{1}{2}-\frac{1}{k}},
\end{align*}
since $q^\frac{1}{k}=p$. For the terms with $(n,p)=1$, we have the following evaluation of the Kloosterman sums, which can be found in \cite[(12.39)]{iwakow}:
\begin{align}
\label{evall} S(n,\beta ;q) =\begin{cases}  2 \big(\frac{\ell}{p} \big)^k q^{\half} {\rm Re} \ \varepsilon_q e\big(\frac{2\ell}{q} \big) &\text{ if } \big(\frac{n\beta}{p}\big)=1, \\
0 &\text{ if } \big(\frac{n\beta }{p}\big)=-1,
\end{cases}
\end{align}
where $\ell^2 \equiv n\beta \bmod q$, $\big( \frac{\ell}{p} \big)$ is the Legendre symbol, and $\varepsilon_q$ equals 1 if $q\equiv 1 \bmod 4$ and $i$ if $q\equiv 3 \bmod 4$. Note that in the first case of (\ref{evall}), the existence of $\ell$ is guaranteed because $n\beta$ is a quadratic residue mod $q$ if and only if it is a quadratic residue mod $p$. Also note that the formula does not depend on the choice of $\ell$. Thus
\begin{align}
&\sum_{\substack{1\le n \le N}} S(n,\beta ;q)  = q^\frac{1}{2} \sum_{\substack{1\le \alpha <p\\ (\frac{\alpha}{p})=1}}  \  \sum_{\substack{1\le n \le N \\ n\beta \equiv \alpha \bmod p}}    \  \sum_{\ell^2 \equiv n\beta \bmod q }   \Big(\frac{\ell}{p} \Big)^k {\rm Re} \  \varepsilon_q e\Big(\frac{2\ell}{q} \Big) + O\big(  N q^{\frac{1}{2}-\frac{1}{k}}  \big).
\end{align}
We will obtain cancellation in only the $n$-sum. Thus to establish (\ref{main2line}), it suffices to prove
\begin{align}
\label{sufff} \sum_{\substack{1\le n \le N \\ n\beta \equiv \alpha \bmod p}}    \  \sum_{\ell^2 \equiv n\beta \bmod q }   \Big(\frac{\ell}{p} \Big)^k {\rm Re} \  \varepsilon_q e\Big(\frac{2\ell}{q} \Big) \ll N q^{-\frac{1}{k}-\delta}
 \end{align}
for some $\delta>0$, uniformly for any integer $1\le \alpha < p$ with $(\frac{\alpha}{p})=1$. 

Let $1\le \gamma \le p$ and $1\le \xi, \omega \le q$ be integers satisfying 
\begin{align*}
&\gamma\equiv \overline{\beta}\alpha \bmod p,\\
&\xi \equiv \overline{ \gamma} \bmod q,\\
&\omega^2 \equiv \beta \gamma \bmod q,
\end{align*}
where overline denotes the multiplicative inverse.
The existence of $\omega$ is guaranteed because $ \beta\gamma \equiv \alpha \bmod p$ is a quadratic residue mod $p$. Since $n\beta\equiv \alpha \bmod p$, we may write 
\begin{align*}
n=  \gamma+ tp
\end{align*}
for some integer $t$. Recall from \cite[Chapter IV.1]{kob} that the power series
\begin{align*}
\sum_{j \ge 0} { \frac{1}{2} \choose j} x^j 
\end{align*}
in ${\mathbb Z}_p[[x]]$ converges in the $p$-adic norm for $x\in p{\mathbb Z}_p$ to a square root of $1+x$. Using this we see that the solutions of the congruence 
\begin{align*}
\ell^2 \equiv n\beta \equiv  \beta\gamma(1+  \xi p t ) \bmod q
\end{align*}
are given by
\begin{align*}
\ell \equiv \pm \omega (1+ c_1 p t + c_2 p^2 t^2 +\ldots +c_{k-1} p^{k-1} t^{k-1}) \bmod q
\end{align*}
for some integers $1\le c_j \le q$ which, for $p$ large enough in terms of $k$, satisfy $(c_j,p)=1$. This also implies $(\frac{\ell}{p}) = (\frac{\pm \omega}{p})$. Thus to prove (\ref{sufff}), it suffices to prove
\begin{align*}
 \sum_{\substack{-\infty \le t \le \infty\\  1\le \gamma + tp \le N}}    E(t)   \ll   Nq^{-\frac{1}{k}-\delta},
\end{align*}
where
\begin{align*}
E(t) &=  e \Big( 2 \omega \Big( \frac{1}{p^k} +  \frac{ c_1 }{p^{k-1}}t +  \frac{ c_2 }{p^{k-2}}t^2 +\ldots + \frac{ c_{k-1} }{p}t^{k-1} \Big)\Big).
\end{align*}
We have 
\begin{align*}
 \sum_{\substack{-\infty \le t \le \infty\\  1\le  \gamma +tp \le N}}   E(t)  =  \sum_{1 \le t \le \frac{N}{p} }    E(t)  + O(1).
\end{align*}
By assumption, $k> \frac{1}{\lambda}$. Suppose also that $\delta\le \lambda k -1$, so that $ Nq^{-\frac{1}{k}-\delta}\ge 1$. Then it suffices to prove
\begin{align}
\label{sufff2} \sum_{ 1 \le t \le \frac{N}{p} }    E(t)   \ll   Nq^{-\frac{1}{k}-\delta}.
\end{align}

Now by Weyl's differencing method, as presented in \cite[Proposition 8.2]{iwakow}, we have that
\begin{align}
\label{lastt} \sum_{ 1 \le t \le \frac{N}{p} }    E(t)   \le   \frac{2N}{p} \Big( \Big(\frac{2N}{p}\Big) ^{-k+1} \sum_{-\frac{N}{p} < j_1, j_2, \ldots, j_{k-2} < \frac{N}{p} } \min\Big\{ \frac{N}{p}, \| \frac{2\omega c_{k-1} (k-1)! j_1j_2\cdots j_{k-2}}{p}\|^{-1} \Big\} \Big)^{2^{2-k}},
\end{align}
where $\| y \|$ is the distance from $y$ to the nearest integer.
Recall that $(\omega c_{k-1},p)=1$. If also $(p,j_2\cdots j_{k-2})=1$ and $p$ is large enough so that  $(p,(k-1)!)=1$, then we have
\begin{align*}
 \sum_{\substack{-\frac{N}{p} < j_1 < \frac{N}{p} \\ (j_1,p)=1} }  \| \frac{2\omega c_{k-1} (k-1)! j_1j_2\cdots j_{k-2}}{p} \|^{-1} \ll \Big(1+\frac{N}{p^2}\Big) p \log p.
\end{align*}
Thus the contribution to (\ref{lastt}) of the terms with $(p,j_1j_2\cdots j_{k-2})=1$ is less than a constant multiple of
\begin{align}
\label{lastt2} \frac{N}{p} \Big( \Big(\frac{N}{p}\Big)^{-k+1} \sum_{-\frac{N}{p} < j_2, \ldots, j_{k-2} < \frac{N}{p} }  \Big(p+\frac{N}{p}\Big) \log p   \Big)^{2^{2-k}} \ll  \frac{N}{p} \Big( \frac{N^2 }{p^3\log p}\Big)^{-2^{2-k}} + \Big(\frac{N}{p}\Big)^{1-2^{2-k}} (\log p)^{2^{2-k}}.
\end{align}
For the right hand side to be $O(Nq^{-\frac{1}{k}-\delta})$ for some $\delta>0$ as required in (\ref{sufff2}), we need $\frac{N^2}{p^3} > p^{\epsilon}$ for some $\epsilon>0$. This is guaranteed by the assumption $k>\frac{3}{2\lambda}$ made in the statement of the theorem. Note that the estimate (\ref{lastt2}) is the very reason for the assumption.

As noted above we have  $\frac{N^2}{p^3} > p^{\epsilon}$, whence $(\frac{N}{p})^{k-2} > p^{\half(k-2)}$. Thus since by assumption $k\ge 4$, we have $(\frac{N}{p})^{k-2} > p$ and so 
\begin{align*}
 \sum_{\substack{-\frac{N}{p} < j_1, j_2, \ldots, j_{k-2} <\frac{N}{p} \\ p|j_1j_2\cdots j_{k-2} } } 1 \ll \frac{1}{p}\Big(\frac{N}{p}\Big)^{k-2}.
\end{align*}
Using this, we see that the contribution to (\ref{lastt}) of the terms with $p|j_1j_2\cdots j_{k-2}$ is less than or equal to
\begin{align*}
\frac{2N}{p}\Big( \Big(\frac{N}{p}\Big)^{-k+1} \sum_{\substack{-\frac{N}{p} < j_1, j_2, \ldots, j_{k-2} <\frac{N}{p} \\ p|j_1j_2\cdots j_{k-2} } } \frac{N}{p} \Big)^{2^{2-k}}\ll N p^{-1-2^{2-k}},
\end{align*}
as required in (\ref{sufff2}).

\section{Proof of Theorem \ref{main}}

We may assume throughout the rest of the paper that $q^{\frac{3}{2}-\eta} < x <q^{\frac{3}{2}+\eta}$, since Theorem \ref{main} is already known for $x \ge q^{\frac{3}{2}+\eta}$. We first reduce Theorem \ref{main} to a problem on estimating a certain sum of Kloosterman sums. This reduction, which is standard, is carried out over the next two lemmas.

\begin{lemma}[Separation of variables] Theorem \ref{main} follows from proving that there exist some fixed positive constants $\eta$ and $\delta$ such that
\begin{align}
\label{lem1} \sum_{uv \equiv a \bmod q} f\Big(\frac{u}{U}\Big)g\Big(\frac{v}{V}\Big)  - \frac{1}{\phi(q)} \sum_{(u,p)=1} \sum_{(v,p)=1} f\Big(\frac{u}{U}\Big)g\Big(\frac{v}{V}\Big) \ll \frac{x^{1-5\delta}}{q}
\end{align} 
for any real numbers $U\ge 1$ and $V\ge 1$ with $1 \le UV \le x$ and any smooth functions $f$ and $g$ compactly supported on $[1,1+x^{-5\delta}]$ with derivatives satisfying
\begin{align}
\label{derbounds}\| f^{(j)} \|_\infty \ll_j x^{5\delta j}, \ \ \ \ \ \| g^{(j)} \|_\infty \ll_j x^{5\delta j}.
\end{align}
\end{lemma}

\proof

We follow \cite{irv}. Opening the divisor function, the left hand side of (\ref{conj}) equals
\begin{align*}
\sum_{\substack{uv\le x\\ uv \equiv a \bmod q}} 1 - \frac{1}{\phi(q)} \sum_{\substack{uv \le x\\ (uv,p)=1}} 1.
\end{align*}
We cover the range of summation by the shorter, almost dyadic, intervals
\begin{align*}
U_i \le u \le (1+ x^{-\frac{3}{2}\delta})U_i, \ \ \ \ V_j\le v \le(1 + x^{-\frac{3}{2}\delta}) V_j
\end{align*}
for $1\le i\le x^{2\delta}$ and $1\le j \le x^{2\delta}$, where $U_i=(1+x^{-\frac{3}{2}\delta})^{i-1}$ and $V_{j}=(1+x^{-\frac{3}{2}\delta})^{j-1}$. Thus Theorem \ref{main} follows from showing
\begin{align}
\label{partitioned} \sum_{\substack{1\le i\le x^{2\delta}\\ 1\le j\le x^{2\delta} \\ U_iV_j \le x}} \Big(  \sum_{\substack{uv\le x\\ uv \equiv a \bmod q \\ U_i \le u \le U_i + x^{-3\delta/2}U_i\\ V_j \le v \le V_j + x^{-3\delta/2}V_j}} 1 - \frac{1}{\phi(q)} \sum_{\substack{uv \le x\\ (uv,p)=1 \\ U_i \le u \le U_i + x^{-3\delta/2}U_i\\ V_j \le v \le V_j + x^{-3\delta/2}V_j}} 1 \Big) \ll \frac{x^{1-\delta}}{q}.
\end{align}
Next we would like to relax the condition $uv\le x$. Note that for any $U_iV_j \le x$, we have that 
\begin{align*}
(U_i + x^{-\frac{3}{2}\delta}U_i)(V_j + x^{-\frac{3}{2}\delta}V_j) < x + 3x^{1-\frac{3}{2}\delta}.
\end{align*}
Thus the left hand side of (\ref{partitioned}) equals
\begin{align}
\label{relax}& \sum_{\substack{1\le i\le x^{2\delta}\\ 1\le j\le x^{2\delta} \\ U_iV_j \le x}}   \Big(  \sum_{\substack{ uv \equiv a \bmod q \\ U_i \le u \le U_i + x^{-3\delta/2}U_i\\ V_j \le v \le V_j + x^{-3\delta/2}V_j}} 1 - \frac{1}{\phi(q)} \sum_{\substack{\\ (uv,p)=1 \\ U_i \le u \le U_i + x^{-3\delta/2}U_i \\ V_j \le v \le V_j + x^{-3\delta/2}V_j}} 1 \Big)
\\ \nonumber &+O\Big( \sum_{\substack{uv \equiv a \bmod q \\ x<uv<x + 3x^{1-3\delta/2}}} 1 + \frac{1}{q}\sum_{\substack{x<uv<x + 3x^{1-3\delta/2}}} 1\Big).
\end{align}
If $\delta$ and $\eta$ are small enough then $x^{1-\frac{3}{2}\delta}\gg q$ and the error term of (\ref{relax}) is $O(x^{1-\delta}/q)$.
Thus to prove Theorem \ref{main}, it suffices to consider only the main term of (\ref{relax}). To this end, 
the indicator function of the interval $U_i \le u \le U_i + x^{-\frac{3}{2}\delta}U_i$ can be approximated by a smooth bump function $f_i(\frac{u}{U_i})$ which equals 0 on the complement of this interval and 1 on the interval $[(1+x^{-5\delta})U_i, (1+x^{-\frac{3}{2}\delta}-x^{-5\delta})U_i]$, and satisfies (\ref{derbounds}). Similarly, the indicator function of the interval $V_j \le v \le V_j + x^{-\frac{3}{2}\delta}V_j$ can be approximated by a smooth bump function $g_j(\frac{v}{V_j})$. For each pair $U_i$ and $V_j$ such that $U_iV_j\le x$, the number of products $uv$ for which $0<f(\frac{u}{U_i})<1$ or $0<g(\frac{v}{V_j})<1$ is $O(x^{1-5\delta})$. Thus the main term of (\ref{relax}) equals
\begin{align*}
\sum_{\substack{1\le i\le x^{2\delta}\\ 1\le j\le x^{2\delta} \\ U_iV_j \le x}}  \Big( \sum_{uv \equiv a \bmod q} f_i\Big(\frac{u}{U_i}\Big)g_j\Big(\frac{v}{V_j}\Big)  - \frac{1}{\phi(q)} \sum_{(u,p)=1} \sum_{(v,p)=1} f_i\Big(\frac{u}{U_i}\Big)g_j\Big(\frac{v}{V_j}\Big) \Big)+ O\Big(   \sum_{\substack{1\le i\le x^{2\delta}\\ 1\le j\le x^{2\delta} \\ U_iV_j \le x}}  \frac{x^{1-5\delta}}{q} \Big),
\end{align*}
provided $\delta$ and $\eta$ are small enough so that $x^{1-5\delta}\gg q$. The error term above is $O(x^{1-\delta}/q)$. Thus to prove Theorem \ref{main}, it suffices to show that
\begin{align}
\label{free} \sum_{uv \equiv a \bmod q} f\Big(\frac{u}{U}\Big)g\Big(\frac{v}{V}\Big)  - \frac{1}{\phi(q)} \sum_{(u,p)=1} \sum_{(v,p)=1} f\Big(\frac{u}{U}\Big)g\Big(\frac{v}{V}\Big) \ll \frac{x^{1-5\delta}}{q}
\end{align} 
for any $U,V$ and any $f,g$ as in the statement of the lemma.
\endproof

Let
\begin{align}
\label{fourier} \hat{f}(\xi) = \int_{-\infty}^{\infty}  f(y) e(-y\xi) \ dy
\end{align}
denote the Fourier transform of $f$, where $e(y)=\exp(2\pi i y)$.

\begin{lemma}[Poisson summation] Theorem \ref{main} follows from proving that there exist some fixed positive constants $\eta$ and $\delta$, depending on $k$, such that
\begin{align}
\label{lem2} \frac{UV}{q^2}  \sum_{\substack{-\infty< n,m < \infty}}   \hat{f}\Big(\frac{nU}{q}\Big) \hat{g}\Big(\frac{mV}{q}\Big) S(a,nm;q) \ll \frac{x^{1-5\delta}}{q}
\end{align}
for any real numbers $U\ge 1$ and $V\ge 1$ with $1 \le UV \le x$ and any smooth functions $f$ and $g$ compactly supported on $[1,1+x^{-5\delta}]$ with derivatives satisfying (\ref{derbounds}).
\end{lemma}

\noindent {\bf Remark.} Weil's bound $S(a,nm;q)\ll q^{\frac{1}{2}}$ and the rapid decay of $\hat{f}$ give an estimate for (\ref{lem2}) which is barely insufficient.

\proof 
It is enough to show that (\ref{lem2}) implies (\ref{lem1}). Using additive characters to pick out the residue class $a \bmod q$, we see that we need to prove that
\begin{align}
\label{free} \frac{1}{q} \sum_{u} \sum_{v}  f\Big(\frac{u}{U}\Big)g\Big(\frac{v}{V}\Big) \sum_{1\le h \le q} e\Big(\frac{h(uv-a)}{q}\Big) - \frac{1}{\phi(q)}  \sum_{(u,p)=1} \sum_{(v,p)=1}  f\Big(\frac{u}{U}\Big)g\Big(\frac{v}{V}\Big)  \ll \frac{x^{1-5\delta}}{q}.
\end{align}
Reordering the $h$-sum above by the greatest common divisor of $h$ and $q$, we write
\begin{align}
\label{mobinv} \sum_{1\le h \le q} e\Big(\frac{h(uv-a)}{q}\Big) = \sum_{0\le r \le k} \ \sum_{\substack{1\le h \le q\\ p^r\| h}} e\Big(\frac{h(uv-a)}{q}\Big) = \sum_{0\le r \le k} \ \summ_{b \bmod p^{k-r}} e\Big(\frac{b(uv-a)}{p^{k-r}}\Big),
\end{align}
where $\summ$ means that summation is restricted to the primitive residue classes and $p^r\| h$ means that $p^r | h$ and $p^{r+1}\nmid h$. Substituting (\ref{mobinv}) into (\ref{free}) and writing the term corresponding to $r=k$ separately, we have that the left hand side of (\ref{free}) equals
\begin{align}
\label{free2} & \frac{1}{q} \sum_{u} \sum_{v} f\Big(\frac{u}{U}\Big)g\Big(\frac{v}{V}\Big) \sum_{0\le r < k} \ \summ_{b \bmod p^{k-r}} e\Big(\frac{b(uv-a)}{p^{k-r}}\Big)\\
\label{free3} &+ \frac{1}{q} \sum_{u} \sum_{v} f\Big(\frac{u}{U}\Big)g\Big(\frac{v}{V}\Big) - \frac{1}{\phi(q)}  \sum_{(u,p)=1} \sum_{(v,p)=1} f\Big(\frac{u}{U}\Big)g\Big(\frac{v}{V}\Big).
\end{align}
Writing $\phi(q)=p^{k-1}(p-1)$, we have that (\ref{free3}) equals
\begin{align}
\label{free3bound} \frac{1}{p^k}\Big( \sum_{u} \sum_{v} -  \frac{1}{1-\frac{1}{p}} \Big( \sum_{u}  - \sum_{p|u} \Big)\Big( \sum_{v}- \sum_{p|v} \Big)\Big) f\Big(\frac{u}{U}\Big)g\Big(\frac{v}{V}\Big) 
\ll \frac{UV}{p^{k+1}}\ll \frac{x}{q^{1+1/k}}.
\end{align}
If $\delta$ is small enough then we have that (\ref{free3bound}) is $O(x^{1-5\delta}/q)$. Thus to establish (\ref{lem1}), it suffices to show that (\ref{free2}) is $O(x^{1-5\delta}/q)$.

By separating $u$ and $v$ into residue classes modulo $p^{k-r}$ and applying Poisson summation, we have for $(b,p)=1$ that 
\begin{align}
\label{pois} &\sum_{u}\sum_{v} f\Big(\frac{u}{U}\Big)g\Big(\frac{v}{V}\Big) e\Big(\frac{buv}{p^{k-r}}\Big)\\
\nonumber & = \frac{UV}{p^{2(k-r)}}  \sum_{s,t \bmod p^{k-r}}  e\Big(\frac{bst}{p^{k-r}}\Big) \sum_{-\infty< n,m < \infty} e\Big(\frac{sn+tm}{p^{k-r}}\Big)  \hat{f}\Big(\frac{nU}{p^{k-r}}\Big) \hat{g}\Big(\frac{mV}{p^{k-r}}\Big).
\end{align}
By (\ref{derbounds}) and integration by parts in (\ref{fourier}), we have that the right hand side of (\ref{pois}) can be restricted, up to an error of $O(q^{-100})$ say, to
\begin{align}
\label{range} |n|< \frac{q^{1+6\delta}}{p^rU}, \ \ \ \ \ |m|< \frac{q^{1+6\delta}}{p^rV}.
\end{align}
By evaluating the $s$-sum on the right hand side of (\ref{pois}), we have that
\begin{align*}
\sum_{u}\sum_{v} f\Big(\frac{u}{U}\Big)g\Big(\frac{v}{V}\Big) e\Big(\frac{buv}{p^{k-r}}\Big)=  \frac{UV}{p^{k-r}}  \sum_{-\infty< n,m < \infty} e\Big(\frac{-nm\overline{b}}{p^{k-r}}\Big)  \hat{f}\Big(\frac{nU}{p^{k-r}}\Big) \hat{g}\Big(\frac{mV}{p^{k-r}}\Big).
\end{align*}
Thus (\ref{free2}) equals
\begin{align}
\label{pois2}  \sum_{0\le r < k} \frac{UV}{p^{2k-r}}  \sum_{-\infty< n,m < \infty}   \hat{f}\Big(\frac{nU}{p^{k-r}}\Big) \hat{g}\Big(\frac{mV}{p^{k-r}}\Big) S(a,nm;p^{k-r}).
\end{align}
By Weil's bound and (\ref{range}), we have that the contribution to (\ref{pois2}) of the terms with $r>0$ is less than a constant multiple of
\begin{align}
\label{triv} \sum_{0< r < k} \frac{UV}{p^{2k-r}}  \sum_{ |n| < \frac{q^{1+6\delta}}{p^rU}} \  \sum_{ |m| < \frac{q^{1+6\delta}}{p^rV}} p^{\frac{k-r}{2}} \ll \sum_{0< r < k} p^{\frac{k}{2}-\frac{3r}{2}+12\delta k} \ll p^{\frac{k}{2}-\frac{3}{2}+12 \delta k}.
\end{align}
This is $O(x^{1-5\delta}/q)$ if $\eta$ and $\delta$ are taken to be small enough. Thus it suffices to prove the bound $O(x^{1-5\delta}/q)$ for only the term with $r=0$ on the right hand side of (\ref{pois2}). This is precisely what was needed to be proved. 
\endproof

\smallskip 

We now use Theorem \ref{main2} to prove Theorem \ref{main}. The goal is to establish (\ref{lem2}). Since $UV\le x$, we suppose by symmetry that $U\le x^\half$. First observe by Weil's bound, $\| \hat{f} \|_\infty < 1$ and (\ref{range}), with $r=0$, that 
\begin{align}
\label{pdiv} \frac{UV}{q^2}  \sum_{\substack{-\infty< n,m < \infty \\ p|nm}}  \Big| \hat{f}\Big(\frac{nU}{q}\Big) \hat{g}\Big(\frac{mV}{q}\Big) S(a,nm;q) \Big| \ll  \frac{UV}{q^2}  \sum_{ \substack{|nm| < \frac{q^{2+12\delta}}{UV} \\ p|nm}}  q^{\frac{1}{2}} \ll  p^{\frac{k}{2}-1+12 \delta k}\log p.
\end{align}
The last bound above uses that $UV\le x <q^{\frac{3}{2}+\eta}$ and $p\le q^\frac{1}{7}$, so that if $\delta$ and $\eta$ are small enough we have $\frac{q^{2+12\delta}}{UV}\gg p$. Thus if $\delta$ and $\eta$ are taken to be small enough, we have that (\ref{pdiv}) is $O(x^{1-5\delta}/q)$. We thus have that the left hand side of (\ref{lem2}) is
\begin{align*}
\frac{UV}{q^2}  \sum_{\substack{-\infty< n,m < \infty}}   \hat{f}\Big(\frac{nU}{q}\Big) \hat{g}\Big(\frac{mV}{q}\Big) &S(a,nm;q) \\= &\frac{UV}{q^2} \sum_\pm \  \sum_{\substack{-\infty< m < \infty\\(m,p)=1}}   \hat{g}\Big(\frac{mV}{q}\Big)  \sum_{\substack{n\ge 1}} \hat{f}\Big(\frac{nU}{q}\Big)S(n,\pm am;q)  + O\Big(  \frac{x^{1-5\delta}}{q} \Big).
\end{align*}
We will obtain cancellation in only the $n$-sum. Bounding the $m$-sum absolutely and using (\ref{range}) again, we see that to establish (\ref{lem2}), it suffices to prove
\begin{align}
\label{suff} \frac{U}{q^{1-6\delta}}   \sum_{\substack{n\ge 1}}  \hat{f}\Big(\frac{nU}{q}\Big)   S(n,\beta;q)  \ll   \frac{x^{1-5\delta}}{q}
\end{align}
uniformly for any integer $1\le \beta < q$ with $(\beta,p)=1$. 
By partial summation and the rapid decay of $\hat{f}$, we have that (\ref{suff}) is bounded by
\begin{align}
\label{suff3} \frac{U}{q^{1-6\delta}}   \sum_{\substack{1 \le N < \frac{q^{1+7\delta}}{U}} }  \Big|  \hat{f}\Big(\frac{(N+1)U}{q}\Big)  -  \hat{f}\Big(\frac{NU}{q}\Big)  \Big| \Big| \sum_{1 \le n \le N} S(n,\beta;q) \Big| .
\end{align}
Since $\| (\hat{f})' \|_\infty < 1$, the difference of $\hat{f}$ values above is less than $\frac{U}{q}$. Using this and Weil's bound, the contribution of the terms in (\ref{suff3}) with $N<\frac{q^{1-7\delta}}{U}$ is $O(x^{1-5\delta}/q)$ if $\eta$ is taken to be small enough.  Now consider the terms with $N\ge \frac{q^{1-7\delta}}{U}$. For $\delta$ and $\eta$ small enough this implies $N\ge q^{\frac{1}{4}-\frac{1}{100}}$, say. By Theorem \ref{main2} then, we have for $k\ge 7$ that
\begin{align*}
\sum_{1 \le n \le N} S(n,\beta;q) \ll N q^{\frac{1}{2}-30\delta}
\end{align*}
for some $\delta>0$ small enough. Note that Theorem \ref{main2} is where this lower bound on $k$ comes from.

We have shown that for $k\ge 7$, the part of the sum (\ref{suff3}) with $N\geq \frac{q^{1-7\delta}}{U}$ is bounded by
\begin{align*}
 \frac{U^2}{q^{2-6\delta}}   \sum_{\substack{ \frac{q^{1-7\delta}}{U} \le N < \frac{q^{1+7\delta}}{U}} }  \Big| \sum_{1 \le n \le N} S(n,\beta;q) \Big| \ll  \frac{U^2}{q^{2-6\delta}}    \frac{q^{2+14\delta}}{U^2} q^{\half-30\delta}\ll q^{\half-10\delta}.
\end{align*}
This is $O(x^{1-5\delta}/q)$ if $\delta$ and $\eta$ are small enough.

\bigskip

{\bf Acknowledgement.} The author thanks D. Mili\'cevi\'c and the anonymous referee for some valuable comments.

\bibliographystyle{amsplain}
\bibliography{divisor}

\end{document}